%version after

%deleting the parenthesis in the first sentence in Proposition 4.2. 

%version 10042007

\input amstex
\def\dd{\hbox{$\,$\rm d}}

\def\lge{\hbox{$\le\kern -0.23em\ge$}}
\def\linesep{\bigskip\nointerlineskip\moveright .3true in
\vbox{\hrule width 5.2true in height .1pt}\vskip2pt
\nointerlineskip\moveright .3true in
\vbox{\hrule width 5.2true in height .1pt}
\bigskip}

\documentstyle{amsppt}
\loadbold
\NoBlackBoxes

\topmatter
\title Explicit Estimate on Primes between Consecutive Cubes
\endtitle
\author Cheng, Yuanyou Furui
\endauthor
\address 
Department of Mathematics, University of North Carolina at Chapel Hill,
Chapel Hill, NC 27599. 
\endaddress
\email YFC\@math.unc.edu 
\endemail
\thanks
I wish to thank the referee for many helpful comments including 
an improved version for Lemma 3.2. 
\endthanks
\abstract 
We give an explicit form of Ingham's Theorem on primes in the short 
intervals, and show that there is at least one prime between every 
two consecutive cubes $x\sp{3}$ and $(x+1)\sp{3}$ if $\log\log x\ge 
15$. \endabstract\endtopmatter

\document
\head 1. Introduction \endhead\noindent
Studies about certain problems in number theory are often 
connected to those about the distribution of the prime numbers; 
problems about the distribution of primes are among the central 
ones in number theory. One problem concerning the distribution 
of primes is the distribution of primes in certain intervals.  
For example, Bertrand's postulate asserts that there is a number 
$B$ such that, for every $x>1$, there is at least one prime number 
between $x$ and $Bx$.  If the interval $[x, Bx]$ is replaced 
by a ``short interval'' $[x, x+x\sp{\theta}]$, then the problem 
is more difficult. 

In 1930, Hoheisel showed that there is at least one prime 
in the above mentioned ``short interval'' with 
$\theta=1-\frac{1}{33,000}$ for sufficiently large $x$'s,
see[13]. Ingham [15], in 1941, proved that there is at least 
one prime in $[x, x+x\sp{3/5+\epsilon}]$, where $\epsilon$ is an 
arbitrary positive number tending to zero whenever $x$ is tending 
to infinity, for ``sufficiently large'' $x$'s.  This implies that 
there is at least one prime between two consecutive cubes 
if the numbers involved are ``large enough.''  One of the better 
results in this direction, conjectured by using the Riemann
Hypothesis, is that there is at least one prime between 
$[x, x+x\sp{1/2+\epsilon}]$ for ``sufficiently large'' $x$'s.  
The latter has not been proved or disproved; though better 
results than Hoheisel's and Ingham's are available.  For example,
one may see [2, 3, 12, 15, 17, 18, 19, 26, 28]. \par

These kinds of results would have many useful applications 
if they were ``explicit'' (with all constants being determined 
explicitly).  For references in other directions with explicit 
results, one can see [4, 8, 22, 23, 24, 25]. To figure out 
the ``sufficiently large'' $x$'s related to $\theta$ as 
mentioned above, one needs to investigate the proof in 
a ``slightly different'' way.  As a starting step in 
this direction, we study the distribution of primes between 
consecutive cubes. In this article, we give an explicit form of 
Ingham's Theorem; specifically, we show that there is at least 
one prime between consecutive cubes if the numbers involved are 
larger than the cubes of $x\sb{0}$ where $x\sb{0}=\exp(\exp(15))$
and we also set $T\sb{0}=\exp(\exp(18))$ throughout this paper
accordingly. \par

Our main task is to prove the Density Theorem or to estimate the number of zeros in the 
strip $\sigma>\frac{1}{2}$ for the Riemann zeta function, see Theorem 1 in the follows. 
We let $\beta=\Re(\rho)$ and $I\sb{\beta}(u)$ be the unit step function at the point 
$u=\beta$; that is, $I\sb{\beta}(u)=1$ for $0\le u\le\beta$ and $I\sb{\beta}(u)=0$ 
for $\beta<u\le 1$. One defines $N(u,T):=\sum\sb{0\le\Im(\rho)<T}I\sb{\beta}(u)$ 
and $N(T):=N(0,T)$. \par

\proclaim{Theorem 1} Let $\frac{5}{8}\le \sigma<1$ and $T\ge T\sb{0}$. One has 
$$N(\sigma, T) 
\le C\sb{D}T\sp\frac{8(1-\sigma)}{3}\log\sp{5} T,$$ where $C\sb{D}:=453472.54$.
\endproclaim\par

\proclaim{Theorem 2} Let $x\ge x\sb{0}$, $h\ge3x\sp{2/3}$ and $C\sb{D}$ be defined 
in Theorem 1. Then $$\psi(x+h)-\psi(x)\ge\,h(1-\epsilon(x)),$$ 
where $$|\epsilon(x)|:=3192.34\exp\left(-\frac{1}{273.79}
\left(\frac{\log\,x}{\log\log\,x}\right)\sp\frac{1}{3}\right).$$\endproclaim

\proclaim{Theorem 3}Let $x\ge\exp(\exp(45))$ and $h\ge3x\sp\frac{2}{3}$. Then 
$$\pi(x+h)-\pi(x)\ge\,h\left(1-3192.34\exp\left(-\frac{1}{283.79}
\left(\frac{\log\,x}{\log\log\,x}\right)\sp\frac{1}{3}\right)\right).$$\endproclaim

\proclaim{Corollary}Let $x\ge\exp(\exp(15))$. Then there is at least one prime 
between each pair of consecutive cubes $x\sp{3}$ and $(x+1)\sp{3}$. 
\endproclaim

The proof of Theorem 1 is delayed until Section 5. We shall prove Theorem 2 and 3 
in Section 2. The proof of Theorem 2 is based on Theorem 1 and Laudau's approximate 
formula, which is in Section 6. Then, it is not difficult to prove Theorem 3 from 
Theorem 2, as shown in Section 2.\par

\head 2. Proof of Theorem 2 and 3\endhead\noindent
From [25], one has 
$$N(T) \le\tfrac{T}{2\pi} \log\,\tfrac{T}{2\pi} 
-\tfrac{T}{2\pi} +\tfrac{7}{8} +0.137\log T +0.443\log\log T +1.588.$$
The following proposition follows straightforward.\par

\proclaim{Proposition 2.1} For $T\ge\,6$, one has $N(T)\le \frac{T\log\,T}{2\pi}$.
\endproclaim

\proclaim{Proposition 2.2}Let $C\sb{D}$ be defined in Theorem 1. Assume that 
the Riemann zeta-function does not vanish for $\sigma>1-z(t)$. Suppose that
$T\sb{0}\le\,T<x\sp{3/8}$. For any $h>0$, one has 
$$\left|\sum\sb{|\Im(\rho)|\le T}\frac{(x+h)\sp{\rho}-x\sp{\rho}}{\rho}\right|
\le\frac{2C\sb{D}T\sp{\frac{8}{3}z(t)}\log x\log\sp{5}T}
{x\sp{z(t)}(\log\,x-\frac{8}{3}\log\,T)}\;h.$$\endproclaim

\noindent{\it Proof}. Notes that
$$\left|\frac{(x+h)\sp{\rho}-x\sp{\rho}}{\rho}\right|
=\left|\int\sb{x}\sp{x+h}u\sp{\rho-1}\dd\,u\right|\le\,hx\sp{\beta-1},$$
where $\beta=\Re(\rho)$ is the real part of $\rho$;
$$x\sp{\beta}=1+\log\,x\int\sb{0}\sp{\beta}x\sp{u}\dd\,u;$$
and
$$\int_{0}^{\beta} x^{u} \dd u =\int_{0}^{1} x^{u} I_{\beta} (u) \dd u,$$
where $I_{\beta}(u)$ is the unit step function or $I_{\beta}(u)=1$
for $0\leq u\leq \beta$ and $I_{\beta}(u)=0$ for $\beta<u\leq 1$.
After interchanging the summation and integration, one has
$$\aligned &\left|\sum\sb{|\Im(\rho)|\le\,T} \frac{(x+h)\sp{\rho}
    	-x\sp{\rho}}{\rho}\right|
	\le\frac{h}{x}\sum\sb{|\Im(\rho)|\le\,T}x\sp{\beta}\\
&\hskip 2 true cm \le\frac{h}{x} \left(\sum\sb{|\Im(\rho)|\le\,T}1
	+\log x\int\sb{0}\sp{1}x\sp{u}
	\left(\sum\sb{|\Im(\rho)|\le\,T} I\sb{\beta}(u)\right)\dd u\right).\\ 
\endaligned\tag{2.1}$$
If the Riemann zeta-function does not vanish in the region $\sigma>1-z(t)$, then 
the expression in the outmost parenthesis in (2.1) is bounded by
$$2N(0,T)+2N(0,T)\log x\int\sb{0}\sp{\frac{5}{8}} x\sp{u}\dd u
   +2\log x\int\sb{\frac{5}{8}}\sp{1-z(t)}x\sp{u}N(u,T)\dd u.\tag{2.2}$$ 
Since $T\ge6$, one can apply Proposition 2.1. The sum of the first two terms in (2.2)
$$2\left(1+\log x\int\sb{0}\sp{5/8} x\sp{u}\dd u\right) N(0,T)
=2x\sp{5/8}N(0,T)\le\frac{x\sp{5/8} T\log T}{\pi}.\tag{2.3}$$
From Theorem 1, one sees that the last term in (2.2) is bounded by 
$$\multline
2C\sb{D}\log\,x\log\sp{5}T\int\sb{5/8}\sp{1-z(t)}x\sp{u}T\sp{\frac{8}{3}(1-\sigma)}\dd\,u\\
=2C\sb{D}T\sp{\frac{8}{3}}\log x\log\sp{5}T\int\sb{5/8}\sp{1-z(t)}
	\left(\frac{x}{T\sp{8/3}}\right)\sp{u}\dd\,u\\
=\frac{2C\sb{D}T\sp\frac{8}{3}\log x\log\sp{5}T}{\log\,x-\frac{8}{3}\log\,T}
	\left(\left(\frac{x}{T\sp\frac{8}{3}}\right)\sp{1-z(t)}
	-\left(\frac{x}{T\sp\frac{8}{3}}\right)\sp\frac{5}{8}\right)\\
=\frac{2C\sb{D}xT\sp{\frac{8}{3}z(t)}\log x\log\sp{5}T}{x\sp{z(t)}
	(\log\,x-\frac{8}{3}\log\,T)}
	-\frac{2C\sb{D}x\sp\frac{5}{8}T\log x\log\sp{5}T}{\log\,x-\frac{8}{3}\log\,T}.\\
\endmultline\tag{2.4}$$
One sees that the sum of the upper bound in (2.3) and the second term on the right side 
in (2.4) is negative. Finally, one combines (2.1) and the first term in the last expression
in (2.4) to finish the proof of Lemma 2.1. \par\bigskip

\noindent{\it Proof of Theorem 2}. From Lemma 9.1 and Proposition 2.2, one sees that 
$$\psi(x+h)-\psi(x)=h+h\epsilon(x),$$ with 
$$\eqalign{
|\epsilon(x)|&\le\frac{1}{h}\left(\left|\sum\sb{|\Im(\rho)|\le\,T\sb{u}}
	\frac{(x+h)\sp{\rho}-x\sp{\rho}}{\rho}\right|+|E(x+h)|+|E(x)|\right) \cr
&\le\, \frac{2C\sb{D}T\sp{\frac{8}{3}z(t)}\log\,x\log\sp{5}T}
	{x\sp{z(t)}(\log\,x-\frac{8}{3}\log\,T)} 
		+10.52\,\frac{(x+h)\log\sp{2}(x+h)}{hT} + \cr
&\qquad +66.976\frac{(x+h)\log\sp{2}T}{hT\log x} + 6\frac{\log\sp{2}T}{hx}. \cr
}\tag{2.5}$$

Let $3x\sp\frac{2}{3}\le\,h$. Also, let 
$$T=T(x):=x\sp\frac{1}{3}\exp\left(\frac{1}{256.59}
	\left(\frac{\log\,x}{\log\log\,x}\right)\sp\frac{1}{3}\right),$$ 
with some undetermined constant $u>1$. 
Then, 
$$\log\,T=\frac{1}{3}\log\,x+\frac{1}{256.59}
	\left(\frac{\log\,x}{\log\log\,x}\right)\sp\frac{1}{3}
\le0.34\log\,x,$$
$$\log\log\,T\le\log\log\,x,$$
and
$$T\sp\frac{8}{3}=x\sp\frac{8}{9}\exp\left(\frac{8}{3\times256.59}
\left(\frac{\log\,x}{\log\log\,x}\right)\sp\frac{1}{3}\right).$$\par

From [10], it is known that the Riemann zeta function does not vanish for $T\ge 1-z(T)$
with $$z(T)=\frac{1}{58.51\log\sp{2/3}T(\log\log T)\sp{1/3}}.$$ Let $Z(x):=z(T(x))$.
Then $$Z(x)\ge\frac{1}{28.51\log\sp{2/3}x(\log\log\,x)\sp{2/3}},$$
$$\multline\left(\frac{x}{T\sp\frac{8}{3}}\right)\sp{z(T)}
	\ge\left(\frac{x\sp\frac{1}{9}}
	{\exp\left(\frac{8}{3\times256.59}
	\left(\frac{\log\,x}{\log\log\,x}\right)\sp\frac{1}{3}\right)}\right)\sp{Z(x)}\\
=\exp\left(\frac{1}{256.59}\left(\frac{\log\,x}{\log\log\,x}\right)\sp\frac{1}{3}
	-\frac{8}{3\times256.59\times28.51\log\sp\frac{1}{3}x(\log\log\,x)\sp\frac{1 or 2}{3}}
	\right).\\
\endmultline$$
It follows that for $x\ge\exp(\exp(45))$ the right side in (2.5) is bounded from above by 
$$\eqalign{
& 3192.34\exp\left(-\frac{1}{273.79}\left(\frac{\log\,x}{\log\log\,x}\right)\sp\frac{1}{3}\right)  
	+1.76\exp\left(-\frac{1}{25s6.6}\left(\frac{\log\,x}{\log\log\,x}\right)\sp\frac{1}{3}\right) + \cr
&\quad +2.59\exp\left(-\frac{1}{256.6}\left(\frac{\log\,x}{\log\log\,x}\right)\sp\frac{1}{3}\right) 
	+\frac{0.24\log\sp{2}+4.65\log\,x+260.48}{x}\,. \cr
}$$
Conclude that one has proved Theorem 2. 

\noindent{\it Proof of Theorem 3}. By definition of $\pi(x)$ and $\psi(x)$, one has
$$\aligned &\pi(x+h)-\pi(x) =\sum\sb{x<p\le x+h} 1\ge\sum\sb{x<p\le x+h} 
	\frac{\log p}{\log x}\\
&=\frac{\psi(x+h)-\psi(x)}{\log x}\ge\frac{h}{\log x} 
	\left(1-\epsilon(x)\right).\\ \endaligned$$
This finishes the proof of Theorem 3.  

\noindent{\it Proof of Corollary}. Let $X=x\sp{3}$ and $h=(x+1)\sp{3}-x\sp{3}$. 
Then $h\ge3x\sp{2}=3X\sp\frac{2}{3}$. By Theorem 3, 
$$\pi\left((x+1)\sp{3}\right)-\pi\left(x\sp{3}\right)\ge\frac{3x\sp{2}}{3\log\,x}
\left(1-\epsilon(x\sp{3})\right)>1.$$
This proves the corollary. 

\head 3. Three auxiliary functions\endhead\noindent
Three auxiliary functions $U\sb{A}$, $V\sb{A}$, and $W\sb{A}$ are introduced 
in this section. For~references, one may see [4], [17], [26].\par\bigskip

\noindent{\bf Definitions of Three Auxiliary Functions.} Let $A$ be a positive integer. 
Define $$U\sb{A}(s)=\sum\sb{n=1}\sp{A} 
\frac{\mu(n)}{n\sp{s}}.$$ Here $\mu$ is the M\"obius $\mu$-function. 
Then, $$V\sb{A}(s)=\zeta(s)U\sb{A}(s)-1,\quad W\sb{A}(s)=1-V\sb{A}\sp{2}(s).$$\par

\proclaim{Lemma 3.1} Let $\nu(n)=\sum\sb{m\le\,A:m|n}\mu(m)$. Then $|\nu(n)|
\le\,d(n)$ and $$V\sb{A}(s)=\sum\sb{n>A}\frac{\nu(n)}{n\sp{s}}.$$
Every non-trivial zero of $\zeta(s)$ is a zero of $W\sb{A}(s)$.\endproclaim

\proclaim{Lemma 3.2} One has $$|V\sb{A}(2+it)|\sp{2}\le\frac{7.9}{A}.$$ 
If $A\ge8$, then both $\Re(W\sb{A}(2+it))$ and $W\sb{A}(2+it)$ does not vanish; 
if $A\ge16$, then $|V\sb{A}(2+it)|\sp{2}<\frac{1}{2}$ and $|W\sb{A}(2+it)|>
\frac{1}{2}$.
\endproclaim

\proclaim{Lemma 3.3} Let $b\sb{1}=5.134$. For $\sigma\ge\frac{1}{4}$ and $t\ge3.297$, 
one has $$\left|V\sb{A}(s)\right|\ll\,t\sp{3/2},$$ and 
$$|W\sb{A}(s)|\le\left(\frac{16}{9}A\sp{3/4}t\sp{3/2}+b\sb{1}A\sp{3/4}t\sp{1/2}\right)
\left(\frac{16}{9}A\sp{3/4}t\sp{3/2}+b\sb{1}A\sp{3/4}t\sp{1/2}+2\right).$$
\endproclaim

\linesep

\noindent{\it Proof of Lemma 3.1}. Easy.  

\noindent{\it Proof of Lemma 3.2}. We may assume $A\ge 5$. Observe
$$\eqalign{ \sum\sb{n>A} \frac{d(n)}{n\sp{2}} 
&= \zeta(2) \sum\sb{m>A} \frac{1}{m\sp{2}} +\sum\sb{m\le A} \frac{1}{m\sp{2}} 
	\sum\sb{n>\frac{A}{m}}\frac{1}{n\sp{2}} \cr
&\le \frac{\zeta(2)}{A} +\sum\sb{m\le A} \frac{1}{m\sp{2}}\frac{1}{\left\lfloor\frac{A}{m}\right\rfloor} \cr
&\le \frac{\zeta(2)}{A} +\sum\sb{m\le \frac{A}{2}+1} \frac{1}{m (A-m)} +\sum\sb{\frac{A}{2} +1<m\le A}\frac{1}{m\sp{2}} \cr
&\le \frac{\zeta(2)}{A} +\frac{1}{A-1} +\frac{4}{A\sp{2}-4}+\frac{\log(A-1)}{A} +\frac{1}{A} \cr
&<\frac{2.8}{A} \cr
}$$
for $A\ge 5$.  

One needs the following proposition.\par

\proclaim{Proposition 3.1}  Let $\sigma\ge\frac{1}{4}$ and $t\ge 3.297$. Then
$$\zeta(s) =\sum\sb{n=1}\sp{\lfloor t\sp{2}\rfloor}\frac{1}{n\sp{s}}+B(s),$$ 
where $|B(s)|\le\,b\sb{1} t\sp{1/2}$ with $b\sb{1}:=5.134$. \endproclaim 

\noindent{\it Proof}. Note that in [4], [17], [26]
$$\zeta(s)=\sum\sb{n=1}\sp{N}\frac{1}{n\sp{s}} 
-s\int\sb{N}\sp{\infty}\frac{u-[u]}{u\sp{s+1}}\dd u+\frac{1}{(s-1)N\sp{s-1}}
\quad (\sigma>0, s\not=1). $$ 
From this, we have
$$\left|\zeta(s)-\sum\sb{n\le N} \frac{1}{n\sp{s}} \right|
\le \frac{N\sp{1-\sigma}}{t}+\sqrt{\frac{1}{\sigma\sb{0}\sp{2}}
+\frac{1}{t\sb{0}\sp{2}}}\frac{t}{N\sp{\sigma}}\quad (\sigma>0, s\not=1).$$
Using this identity, Proposition 3.1 follows.  

\noindent{\it Proof of Lemma 3.3}. Using its definition, one sees that 
$$|U\sb{A}(s)|\le\sum\sb{n=1}\sp{A}\frac{1}{n\sp{\sigma}}.$$ 
If $0<\sigma<1$, one gets 
$$|U\sb{A}(s)|\le\int\sb{0}\sp{A}\frac{1}{u\sp{\sigma}}\dd\,u 
=\frac{A\sp{1-\sigma}}{1-\sigma};$$ if $\sigma\ge\,1$, one has 
$$|U\sb{A}(s)|\le\sum\sb{n=1}\sp{A}\frac{1}{n}\le\log\,A+1.$$ 
For $\sigma\ge\frac{1}{4}$, one obtains 
$$|U\sb{A}(s)|\le\max\left\{\frac{4}{3}A\sp{3/4}, \log\,A+1\right\}
\le\frac{4}{3}A\sp{3/4}.\tag{3.2}$$\par

Similarly, one gets
$$\sum\sb{n=1}\sp{\lfloor\,t\sp{2}\rfloor}\frac{1}{n\sp{\sigma}}
\le\frac{4}{3} t\sp{3/2}.$$  Combining with the result in Proposition
3.2, one has $$|\zeta(s)|\le\frac{4}{3}t\sp{3/2}+b\sb{1}t\sp{1/2},\tag{3.3}$$
for $\sigma\ge\frac{1}{4}$ and $t\ge 3.297$.\par

Recalling the definition of $V\sb{A}(s)$ one has 
$$|V\sb{A}(s)|\le|\zeta(s)||U\sb{A}(s)|+1;$$ 
from (3.1), one gets 
$$|W\sb{A}(s)|\le|\zeta(s)|\,|U\sb{A}(s)|(2+|\zeta(s)||U\sb{A}(s)|).$$ 
Conclude with (3.2) and (3.3) that one proves Lemma 3.3. \bigskip
 
\head 4. Representing the number of zeros by an integral\endhead \noindent
\noindent{\bf Notation $N\sb{F}(\sigma, T)$.} Let $F(s)$ be a complex function
and $T>0$. The notation $N\sb{F}(\sigma, T)$ expresses the number of zeros 
in the form $\beta+i\gamma$ for $F(s)$ with $\sigma\le\beta$ and $0\le\gamma<T$.\par

It is well known that $\zeta(s)$ does not vanish for $\sigma\ge 1$; so one may 
restrict our discussion to $\sigma<1$.\par

\proclaim{Lemma 4.1} Let $T\sb{1}=14$ and $A\ge16$. Then for $\sigma\sb{0}<\sigma<1$ 
and $T\ge\,T\sb{1}$, one has 
$$N\sb{\zeta}(\sigma; T)\le\frac{1}{\sigma-\sigma\sb{0}}
\left(\frac{1}{2\pi}\int\sb{T\sb{1}}\sp{T}\left|V\sb{A}
\left(\sigma\sb{0}+it\right)\right|\sp{2}\dd t+\frac{(594?)16T}{2\pi\,A}+1+c\sb{A}(T)\right),$$ 
with $$\multline
c\sb{A}(T)
:=\frac{\log\left(\frac{16}{9}A\sp{3/4}\Big(T+\frac{7}{4}\Big)\sp{3/2}
	+b\sb{1}A\sp{3/4}\Big(T+\frac{7}{4}\Big)\sp{1/2}\right)}{\log(7/6)}\\
+\frac{\log\left(\frac{16}{9}A\sp{3/4}\Big(T+\frac{7}{4}\Big)\sp{3/2}
	+b\sb{1}A\sp{3/4}\Big(T+\frac{7}{4}\Big)\sp{1/2}+2\right)}{\log(7/6)}
	+\frac{\log\,2}{\log(7/6)}.\\
\endmultline$$
\endproclaim

\proclaim{Corollary}Let $A\le\frac{595}{594}T$ and $T\ge\exp(\exp(18))$. Then 
$$c\sb{A}(T)\le\,29.193\log\,T+11.978.$$
\endproclaim
\linesep

\noindent{\bf Notation $N\sb{F}(\sigma; T, T\sb{1})$.} Let $F(s)$, $\sigma$, 
and $T$ as in the last definition. The notation $N\sb{F}(\sigma; T, T\sb{1})$ 
expresses the number of zeros in the form $\beta+i\gamma$ for $F(s)$ 
with $\sigma\le\beta$ and $T\sb{1}\le\gamma<T$.\par\medskip

Be definition, one sees that $N\sb{F}(\sigma;T,T\sb{1})=N\sb{F}(\sigma,T)
-N\sb{F}(\sigma,T\sb{1})$ for any complex function $F$. Note here, see [9] or 
[17], that there is no zero for the Riemann zeta function $\zeta(\sigma+it)$
for $0\le\,t\le14$. If one takes $T\sb{1}=14$, then $N\sb{\zeta}(\sigma;T,T\sb{1})
=N\sb{\zeta}(\sigma, T)$.\par

For an analytic function, a zero is isolated and the number of zeros in any
compact region is finite. Fix $\sigma$ and $T$. Let $\epsilon\sb{1}$, 
$\epsilon\sb{2}$, and $\epsilon\sb{3}$ be sufficiently small positive numbers 
and $\lambda=\sigma-\epsilon\sb{1}$ and $T\sb{2}=T+\epsilon\sb{2}$. One may 
assume that $\lambda$ is not the real part and $T\sb{2}$ is not the imaginary 
parts of any zeros for the function $W\sb{A}(s)$. Recalling the second part 
of Lemma 3.1, one gets the following proposition.\par

\proclaim{Proposition 4.1} Let $T\sb{1}=14$ and $\epsilon\sb{1}$ and
$\epsilon\sb{2}$ be small positive numbers such that $\lambda=\sigma-\epsilon\sb{1}$ 
is not the real part and $T\sb{2}=T+\epsilon\sb{2}$ is not the imaginary part 
of any zero for the function $W\sb{A}(s)$. Then $$N\sb{\zeta}(\sigma,T)\le\,
N\sb{W\sb{A}}(\lambda;T\sb{2},T\sb{1}).$$ 
\endproclaim

Let $\sigma\sb{0}<\lambda$. Since $N\sb{W\sb{A}}(\lambda; T\sb{2}, T\sb{1})$ 
is a non-increasing function of $\lambda$ by the definition, one sees that
$$N\sb{W\sb{A}}(\lambda; T\sb{2}, T\sb{1})\le\frac{1}{\lambda-\sigma\sb{0}}
\int\sb{\sigma\sb{0}}\sp{\lambda}N\sb{W\sb{A}}(\rho;T\sb{2},T\sb{1})\dd\rho.$$ 
Noting that $$\int\sb{\sigma\sb{0}}\sp{\lambda}N\sb{W\sb{A}}(\rho;T\sb{2},
T\sb{1})\dd\rho\le\int\sb{\sigma\sb{0}}\sp{2}N\sb{W\sb{\bold A}}(\rho;T\sb{2},
T\sb{1})\dd\rho,$$ one has the next proposition.

\proclaim{Proposition 4.2} Let $\sigma\sb{0}<\lambda<1$ and $T\sb{2} >T\sb{1}$.
Assume that $\lambda$ is not the real part and $T\sb{2}$ is not the imaginary part 
of any zero for $W\sb{A}(s)$. Then 
$$N\sb{W\sb{A}}(\lambda; T\sb{2}, T\sb{1})\le\frac{1}{\lambda-\sigma\sb{0}}
\int\sb{\sigma\sb{0}}\sp{2}N\sb{W\sb{A}}(\rho; T\sb{2}, T\sb{1})\dd\rho.$$
\endproclaim

Using the arguments in [26, p.213 and p.220], one gets the following result.\par

\proclaim{Proposition 4.3} Let $\frac{1}{2}<\sigma\sb{0}<2$. Assume that 
$T\sb{2}$ is not the imaginary part of any zero for $W\sb{A}(s)$. Also, 
let ${\Cal N}\sb{k}$ be the number of zeros for $\Re(W\sb{A}(s))$ on the segment 
between $\sigma\sb{0}+it$ and $2+it$ on the line $t=T\sb{k}$ for $k=1$ and $2$ 
respectively. Then 
$$\int\sb{\sigma\sb{0}}\sp{2} N\sb{W\sb{A}}(\rho;T\sb{2},T\sb{1})\dd\rho
\le\frac{1}{2\pi}\int\sb{T\sb{1}}\sp{T\sb{2}}\log\left(\frac{|W\sb{A}(\sigma\sb{0}+it)|}
{|W\sb{A}(2+it)|}\right)\dd t+\frac{{\Cal N}\sb{1}+{\Cal N}\sb{2}}{2}+1.\tag{4.5}$$
\endproclaim

The following result can be found in [4].\par

\proclaim{Proposition 4.4} Suppose that $s\sb{0}$ is a fixed complex number and 
$f$ is a complex function non-vanishing at $s\sb{0}$ and regular for $|s-s\sb{0}|<R$ 
for positive number $R$.  Let $0<r<R$ and $M\sb{f}=\max\sb{|s-s\sb{0}|=R}|f(s)|$.  
Then the number of zeros of $f$ in $|s-s\sb{0}|\le r$, denoted by ${\Cal N}\sb{f}$, 
multiple zeros being counted according to their order of multiplicity satisfies 
the following inequality.  $${\Cal N}\sb{f} 
\le\frac{\log M\sb{f}-\log |f(s\sb{0})|}{\log R-\log r}.$$
\endproclaim

\noindent{\it Proof of Lemma 4.1}. From Proposition 4.1, 4.2, and 4.3, one has
$$N(\sigma, T)\le\frac{1}{\lambda-\sigma\sb{0}}
\left(\frac{1}{2\pi}\int\sb{T\sb{1}}\sp{T\sb{2}}\log\frac{|W\sb{A}(\sigma\sb{0}+it)|}
	{|W\sb{A}(2+it)|}\dd t+\frac{{\Cal N}\sb{1}+{\Cal N}\sb{2}}{2}+1\right),
\tag{4.6}$$ where $\lambda=\sigma-\epsilon\sb{1}$ as $T\sb{2}=T+\epsilon\sb{2}$ 
as in Proposition 4.1.\par

Clearly, we have  
$$\Re(W\sb{A}(\sigma+it))
=\frac{1}{2}\Big(W\sb{A}(\sigma+it)+W\sb{A}(\sigma-it)\Big).$$  
The number of zeros of $\Re(W\sb{A}(s))$ on ${\Cal S}\sb{k}$ are the same 
for~the~following~regular~functions 
$$W\sb{0}\sp{(k)}(s) =\frac{1}{2} \left(W\sb{A}(s+iT\sb{k})
+W\sb{A}(s-iT\sb{k})\right)$$ on the real axis between 
$\sigma\sb{0}$ and $2$.\par

First, one applies Proposition 4.4 to estimate the number ${\Cal N}\sb{k}\sp{\prime}$ 
of zeros for $W\sb{0}\sp{(k)}(s)$ in $|s-2|\le \frac{3}{2}$. It is obvious that 
${\Cal N}\sb{k}\sp{\prime}\le{\Cal N}\sb{k}$. One takes $s\sb{0}=2$, $R=\frac{7}{4}$, 
$r=\frac{3}{2}$. Recalling (4.7) and Lemma 3.3, one acquires 
$$\multline
\max\sb{|s-2|=\frac{3}{2}}\left|W\sb{0}\sp{(k)}(s)\right|
\le\frac{1}{2}(\max\sb{|s-2|=\frac{3}{2}}\left|W\sb{A}(s+iT\sb{k})\right|
	+\max\sb{|s-2|=\frac{3}{2}}\left|W\sb{A}(s-iT\sb{k})\right|)\\
\le\max\sb{|s-2|=\frac{3}{2}}\left|W\sb{A}(s+iT\sb{k})\right|
	\le\,W\sb{A}\sp{(1)}\left(T\sb{k}+\frac{7}{4}\right)
	\le\,W\sb{A}\sp{(1)}\left(T\sb{2}+\frac{7}{4}\right),\\
\endmultline$$ 
where $W\sb{A}\sp{(1)}(t)$ is the upper bound of $|W\sb{A}(s)|$ in Lemma 3.3.
Letting $\epsilon\sb{2}$ tend to zero, one sees that 
$$\max\sb{|s-2|=\frac{3}{2}}|W\sb{0}\sp{(k)}(s)|
	\le\,W\sb{A}\sp{(1)}\left(T+\frac{7}{4}\right).$$
Also, recall that $|W\sb{0}(2+it)|>\frac{1}{2}$ from Lemma 3.2. This implies
$${\Cal N}\sb{k}\le\,c\sb{A}(T):=\frac{\log\,W\sb{A}\sp{(1)}
\left(T+\frac{7}{4}\right)-\log(1/2)}{\log(7/6)},$$  
for $k=1$ and $2$. Hence, $$\frac{{\Cal N}\sb{1}+{\Cal N}\sb{2}}{2}\le
c\sb{A}(T).\tag{4.8}$$\par

Now, transform the integral in (4.6) into a one involving the function 
$V\sb{A}(s)$ instead of $W\sb{A}(s)$. \par

Recalling the definition of $W\sb{A}(s)$, using the triangular inequality 
in the form $|x-y|\le |x|+|y|$, and noting that $\log(1+x)\le x$ for $x>0$,  
one has 
$$\aligned&\log|W\sb{A}(\sigma\sb{0}+it)|
	=\log|1-V\sp{2}(\sigma\sb{0}+it)| \\
&\le \log(1+|V\sb{A}(\sigma\sb{0}+it)|\sp{2})
\le|V\sb{A}(\sigma\sb{0}+it)|\sp{2}.\\ 
\endaligned\tag{4.9}$$\par

Also, by the triangular inequality in the form $|x-y|\ge|x|-|y|$, one gets 
$|1-V\sb{A}\sp{2}(1+it)|\ge 1-|V\sb{A}(1+it)|\sp{2}$.  Using the increasing 
property of the logarithmic function, one sees that $\log|1-V\sb{A}\sp{2}(2+it)|
\ge\log(1-|V\sb{A}(2+it)|\sp{2})$.  It follows that 
$$-\log|W\sb{A}(2+it)|=-\log|1-V\sb{A}\sp{2}(2+it)|
\le-\log(1-|V\sb{A}(2+it)|\sp{2}).$$  
From the last part of Lemma 3.2, one sees that $|V\sb{A}(2+it)|\sp{2}<\frac{1}{2}$
since $A\ge16$. Applying $-\log(1-x)<2x$ for $0<x<\frac{1}{2}$, one acquires 
$$-\log|W\sb{A}(2+it)|\le\,2|V\sb{A}(2+it)|\sp{2}<\frac{7.9}{A}.\tag{4.10}$$\par

Combining (4.9) and (4.10), one obtains 
$$\log\left(\frac{|W\sb{A}(\sigma\sb{0}+it)|}{|W\sb{A}(2+it)|}\right)
\le\left|V\sb{A}\left(\sigma\sb{0}+it\right)\right|\sp{2}+\frac{7.9}{A}.$$\par

Letting $\epsilon\sb{2}$ tend to zeros in (4.6), one gets
$$N(\sigma, T)\le\frac{1}{\lambda-\sigma\sb{0}}
\left(\frac{1}{2\pi}\int\sb{T\sb{1}}\sp{T} 
	\left|V\sb{A}\left(\sigma\sb{0}+it\right)\right|\sp{2}\dd\,t 
	+\frac{7.9T}{2\pi\,A}+1+c\sb{A}(T)\right)\,.$$ 
Finally, letting $\epsilon\sb{1}$ tend to zero, one obtains
$$N\sb{\zeta}(\sigma, T)\le\frac{1}{\sigma-\sigma\sb{0}}
\left(\frac{1}{2\pi}\int\sb{T\sb{1}}\sp{T} 
	\left|V\sb{A}\left(\sigma\sb{0}+it\right)\right|\sp{2}\dd\,t 
	+\frac{7.9T}{2\pi\,A}+1+c\sb{A}(T)\right)\,.$$ 
This proves Lemma 4.1. \bigskip

\head 5. The Proof of Theorem 2\endhead\noindent 
To estimate the integral in Lemma 4.1, one studies the following functions.\par\bigskip

\noindent{\bf Definition of ${\Cal V}\sb{\sigma}(t)$.} Let $t\ge 0$. 
Define $${\Cal V}\sb{\sigma}(t) 
=\int\sb{0}\sp{t}|V\sb{A}(\sigma+iy)|\sp{2}\dd\,y.$$\par

One needs an explicit upper bound for the Riemann zeta function on the line
$\sigma=\frac{1}{2}$, for which we summarize the Corollary and Theorem 1, 2, 
and 3 from [6] into the following lemma. \par

\proclaim{Lemma 5.1}One has 
$$\left|\zeta\left(\frac{1}{2}+it\right)\right|\le\,C\,t\sp{\alpha}
\log\sp{\beta}(t+e)+D$$ for any $t>0$, where $C=3$, $\alpha=\frac{1}{6}$, 
$\beta=1$, and $D=2.657$.
\endproclaim

For $\sigma=\frac{1}{2}$ and $1+\delta$ with the value of $\delta>0$ being 
determined later, one has Lemma 5.2 and 5.3. \par
 
\proclaim{Lemma 5.2} Let $C$, $\alpha$ and $\beta$ as defined in \hbox{\rm Lemma 5.1}, 
$A\ge 16 $, and $0\le t<\infty$. Then, 
$${\Cal V}\sb{1/2}(t)\le D\sb{1} t\sp{2\alpha+1} 
\log\sp{2\beta}(t+e)+D\sb{2} t\sp{2\alpha} \log\sp{2\beta}(t+e)+D\sb{3} t+D\sb{4},$$ 
where $D\sb{1}:=4C\sp{2}(\log\,A+1)$, $D\sb{2}:=16C\sp{2}A(\log\,A+4)$, 
$D\sb{3}:=4D\sp{2}(\log\,A+1)$, and $D\sb{4}:=16D\sp{2}A(\log\,A+4)$.
\endproclaim

%%%  Impossible to combine them together, since $0\le t<\infty$. 

\proclaim{Lemma 5.3}. Let $0<\delta\le 1$,  $A\ge16$, and $0\le t<\infty$. Then 
$${\Cal V}\sb{1+\delta}(t)\le D\sb{5} t+D\sb{6},$$
where $$D\sb{5}:=0.206\frac{\log\sp{3}A+3\log\sp{2}A+6\log\,A+6}{A\sp{1+2\delta}},$$
and $$\aligned
D\sb{6}&:=\frac{0.264(1+\delta)}{A\sp{\delta}}
	\left(\frac{\log\sp{3}A}{\delta}+\frac{3\log\sp{2}A}{\delta\sp{2}}
	+\frac{6\log\,A}{\delta\sp{3}}+\frac{6}{\delta\sp{4}}\right)\\
&\qquad\qquad+\frac{4.012}{A\sp{2\delta}}
	\left(\frac{\log\sp{2}A}{\delta\sp{2}}+\frac{2\log\,A}{\delta\sp{3}}
	+\frac{1}{\delta\sp{4}}\right)\\
&\qquad\qquad\qquad+\frac{16.020(1+\delta)}{A\sp{\delta}}
	\left(\frac{\log\sp{2}A}{\delta}+\frac{2\log\,A}{\delta\sp{2}}
	+\frac{2}{\delta\sp{3}}\right).\\
\endaligned$$
\endproclaim

The proofs of Lemma 5.2 and 5.3 will be given in Section 8. One needs another 
auxiliary function $H(s)$.\par\bigskip

\noindent{\bf Definition of $H(s)$.} Let $\sigma>\frac{1}{2}$ and $t>0$. 
Denote $$H(s):=H\sb{A,\tau}(s):=\frac{s-1}{s\cos\left(\frac{s}{2\tau}\right)} 
V\sb{A}(s), \tag{5.1}$$  
where $V\sb{A}(s)$ is defined in section 3 and $\tau$ is 
a parameter with positive value.\par

The function $H(s)$ has a close relation to the function $V(s)$, as shown in the 
following Lemma 5.4.\par

\proclaim{Lemma 5.4} Let $\frac{1}{2}\le\sigma\le\,2$ and $\tau\ge\,e$. 
Then, for $t>0$, $$|H(s)|<2\,e\sp{-\frac{t}{2\tau}} |V\sb{A}(s)|; $$ 
for $t>14$, 
$$|V\sb{A}(s)|<\,\sqrt{\frac{200}{197}}\,
e\sp{\frac{t}{2\tau}}|H(s)|.$$
\endproclaim\bigskip 

One then uses Lemma 5.2 and 5.3 to give estimates as in Lemma 5.5 on ${\Cal H}(\sigma)$
defined as follows.\par\bigskip

\noindent{\bf Definition of ${\Cal H}(\sigma)$.} Let $\sigma>\frac{1}{2}$. Define 
$${\Cal H}(\sigma)=\int\sb{-\infty}\sp{\infty}|H(s)|\sp{2}\dd t. \tag{5.2}$$\par

\proclaim{Lemma 5.5} For $T\ge T\sb{0}$, one has 
$${\Cal H}\left(\frac{1}{2}\right) \le A\sb{1}T\sp{4/3}\log\sp{3}T, $$
and  
$${\Cal H}\left(1+\frac{\omega}{\log T}\right) \le A\sb{2} \log\sp{4} T, $$ 
where 
$$A\sb{1} =685.026\;\kappa\sp{4/3} +2061.486\kappa\sp{1/3} +0.000001\kappa +0.001, $$ 
and 
$$\eqalign{ A\sb{2} 
&=\frac{144.001}{\pi\sp{2}e\sp{\omega}} 
	\left(\frac{\eta\sp{3}}{\omega}+\frac{3\eta\sp{2}}{\omega\sp{2}}
	+\frac{6\eta\sp{3}}{\omega\sp{3}}+\frac{6\eta\sp{3}}{\omega\sp{4}}\right) 
		+\frac{4}{e\sp{2\omega}} \left(\frac{\eta\sp{2}}{\omega\sp{2}}+\frac{2\eta}{\omega\sp{3}}
			+\frac{1}{\omega\sp{4}}\right) + \cr
&\qquad +\frac{4.689\kappa}{\pi\sp{2}e\sp{2\omega}} +\frac{8.001}{\pi\sp{2}e\sp{\omega}}
	\left(\frac{\eta\sp{2}}{\omega}+\frac{2\eta}{\omega\sp{2}}
	+\frac{2}{\omega\sp{3}}\right), \cr
}$$
with $\eta=1.000001$. 
\endproclaim

One has the following Corollary by taking $\omega=1.598$ and $\kappa=1.501$. The justification of these 
choices of constants will be given in the proof of Theorem 1. \par

\proclaim{Corollary} Let $T\ge T\sb{0}$. Then 
$${\Cal H}\left(\frac{1}{2}\right)\le A\sb{1} T\sp{4/3}\log\sp{3}T\quad\hbox{and}\quad
{\Cal H}\left(1+\frac{1.598}{\log\,T}\right)\le A\sb{2}\log\sp{4}T,$$
with $A\sb{1}=3537.613$ and $A\sb{2}=78.383$. 
\endproclaim

One may transform the estimates on ${\Cal H}(\sigma)$ for $\sigma=\frac{1}{2}$ and 
$1+\delta$ to any $\sigma$ between by the following lemma,  which is due to 
Hardy, Ingham and P\'olya, see [4].\par

\proclaim{Lemma 5.6} If $H(s)$ is regular and bounded for $\sigma\sb{1}\le 
\sigma \le \sigma\sb{2}$, and the integral $${\Cal H}(\sigma) 
=\int\sb{-\infty}\sp{\infty} |H(\sigma+i\,t)|\sp{2}\dd t$$ 
exists, and converges uniformly for $\sigma\sb{1}
\le \sigma\le\sigma\sb{2}$, and 
$$\lim\sb{|t|\to\infty} |H(s)|=0$$ uniformly for $\sigma\sb{1}
\le \sigma\le\sigma\sb{2}$, then for any positive number $T$,  
$${\Cal H}(\sigma)\le 
\{{\Cal H}(\sigma\sb{1})\}\sp\frac{\sigma\sb{2}-\sigma}
     {\sigma\sb{2}-\sigma\sb{1}}
\{ {\Cal H}(\sigma\sb{2})\}\sp\frac{\sigma-\sigma\sb{1}}
{\sigma\sb{2}-\sigma\sb{1}}.$$
\endproclaim 

The proofs of Lemma 5.5 is given in Section 9.\par\bigskip

\noindent{\it Proof of Theorem 1}. Applying Lemma 5.6 to the function 
$H(s)$ with $\sigma\sb{1}=1/2$ and $\sigma\sb{2}=1+\delta$ with any positive $\delta$, 
one obtains $${\Cal H}(\sigma) 
\le\,A\sb{1}\sp\frac{2(1+\delta-\sigma)}{1+2\delta}A\sb{2}\sp\frac{2\sigma-1}{1+2\delta} 
	T\sp\frac{8(1+\delta-\sigma)}{3(1+2\delta)}\log\sp{3+\frac{2\sigma-1}{1+2\delta}}T
\le\,A\sb{1}A\sb{2}T\sp\frac{8(1+\delta-\sigma)}{3}\log\sp{4}T.$$

Now, from Lemma 5.4, one gets 
$|V\sb{A}(s)|\sp{2}\le\frac{200}{197}e\sp{t/\tau}|H(s)|\sp{2}$, or, with $\tau=\kappa T$ 
for $\kappa\ge\frac{e}{T\sb{0}}$, 
$$|V\sb{A}(s)|\sp{2}\le\frac{200}{197}e\sp{1/\kappa}|H(s)|\sp{2}.$$
For the integral in Lemma 4.1, one obtains 
$$\multline\int\sb{T\sb{1}}\sp{T}|V\sb{A}(\sigma\sb{0}+it)|\sp{2}\dd t 
	\le\frac{200}{197}e\sp\frac{1}{\kappa}\int\sb{T\sb{1}}\sp{T}
	|H(\sigma\sb{0}+it)|\sp{2}\dd t\\
\le\frac{200}{197}e\sp\frac{1}{\kappa}\int\sb{0}\sp{\infty}|H(\sigma\sb{0}+it)|\sp{2}\dd t
	=\frac{100}{197} e\sp\frac{1}{\kappa} {\Cal H}(\sigma\sb{0})\\
\le\frac{100}{197} e\sp\frac{1}{\kappa}A\sb{1}A\sb{2}T\sp\frac{8(1+\delta-\sigma\sb{0})}{3} 
	\log\sp{4}T.\\
\endmultline$$
Recalling Lemma 4.1, one sees that 
$$N\sb{\zeta}(\sigma,T) 
\le\frac{100e\sp\frac{1}{\kappa}}{394\pi(\sigma-\sigma\sb{0})}A\sb{1}A\sb{2}
	T\sp\frac{8(1+\delta-\sigma\sb{0})}{3}\log\sp{4}T
	+\frac{1}{\sigma-\sigma\sb{0}}\left(\frac{16 T}{2\pi\,A}+1+c\sb{A}(T)\right).$$

Note that $A\le\big(1+\frac{1}{T\sb{0}}\big) T$ and $\delta=\frac{\omega}{\log T}$ as in the proof of 
Lemma 5.5. Also, let $\sigma\sb{0}=\sigma-\frac{\nu}{\log T}$ for another positive 
constant $\nu$. It follows that 
$$\multline N\sb{\zeta}(\sigma,T) 
\le\frac{100e\sp{\frac{1}{\kappa}+\frac{8(\omega+\nu)}{3}}}{394\pi\nu}
	A\sb{1}A\sb{2}T\sp\frac{8(1-\sigma)}{3}\log\sp{5}T
	+\frac{16\log T}{2\pi\nu}+\frac{\log T}{\nu}
	+\frac{c\sb{A}(T)\log\,T}{\nu}\\
\le\,C\sb{D}T\sp\frac{8(1-\sigma)}{3}\log\sp{5}T,\\
\endmultline$$
with 
$$C\sb{D}:=\frac{100e\sp{\frac{1}{\kappa}+\frac{5\omega}{3}+\frac{8\nu}{3}}}{394\pi\nu}
	A\sb{1}A\sb{2}+\frac{1}{\log\sp{4}T\sb{0}}\left(\frac{16}{2\pi\nu}
	+\frac{1}{\nu}\right)+\frac{c\sb{A}(T)/\log\,T}{\nu\log\sp{3}T\sb{0}}.$$ 
The first term in $C\sb{D}$ is the major one; one may sub-optimize it in order to 
sub-optimize $C\sb{D}$. Note that 
$$\multline
C\sb{D}\approx \frac{100e\sp{\frac{1}{\kappa}+\frac{5\omega}{3}+\frac{8\nu}{3}}}
	{394\pi\nu}(685.026\kappa\sp{4/3} +2061.486\kappa\sp{1/3})\times\\
\left(\frac{144}{\pi\sp{2}e\sp{\omega}}\left(\frac{1}{\omega}+\frac{3}{\omega\sp{2}}
	+\frac{6}{\omega\sp{3}}+\frac{6}{\omega\sp{4}}\right)+
	\frac{4}{e\sp{2\omega}}\left(\frac{1}{\omega\sp{2}}
	+\frac{2}{\omega\sp{3}}+\frac{1}{\omega\sp{4}}\right)\right).\\
\endmultline$$
To optimize the factor $\frac{\exp\frac{8\nu}{3}}{\nu}$, one takes $\nu=\frac{3}{8}$, 
to sub-optimize the factor $$e\sp{1/\kappa}(200.593\kappa\sp{4/3}+603.656\kappa\sp{1/3})$$ 
in $C\sb{D}$, one let $\kappa=1.501$, and to sub-optimize the factor 
$$e\sp{5\omega/3}
\left(\frac{144}{\pi\sp{2}e\sp{\omega}}\left(\frac{1}{\omega}+\frac{3}{\omega\sp{2}}
	+\frac{6}{\omega\sp{3}}+\frac{6}{\omega\sp{4}}\right)+
	\frac{4}{e\sp{2\omega}}\left(\frac{1}{\omega\sp{2}}
	+\frac{2}{\omega\sp{3}}+\frac{1}{\omega\sp{4}}\right)\right),$$
one chooses $\omega=1.598$. With these choices of constants, one gets the Corollary 
of Lemma 5.5. From the Corollary, one justifies the choice of $T\sb{0}$. With computation, 
one finishes the proof of Theorem 1.\par\bigskip

\head 6. Estimates involving the divisor function\endhead
\proclaim{Lemma 6.1} 
Let $\delta>0$ and $\log\log\,N\ge18$. Then
$$\sum\sb{N<n}\frac{d\sp{2}(n)}{n\sp{2+2\delta}}
\le\frac{0.206}{N\sp{1+2\delta}}\left(\log\sp{3}N+3\log\sp{2}N+6\log\,N+6\right).$$
\endproclaim

\proclaim{Lemma 6.2}
Let $\delta>0$ and $\log\log\,N\ge18$. Then
$$\sum\sb{N<n}\sum\sb{N<m<n}\frac{d(m)d(n)}{(mn)\sp{1+\delta}}
\le\frac{1.003}{N\sp{2\delta}}
\left(\frac{\log\sp{2}N}{\delta\sp{2}}+\frac{2\log\,N}{\delta\sp{3}}
+\frac{1}{\delta\sp{4}}\right).$$
\endproclaim

\proclaim{Lemma 6.3}
Let $\delta>0$ and $\log\log\,N\ge18$. Then $$\aligned
&\sum\sb{N<m}\frac{1}{m\sp{1+\delta}}
\sum\sb{N<n<m}\frac{d(m)d(n)}{(mn)\sp{1/2}\log(m/n)}\\
&\qquad\le0.066\frac{1+\delta}{N\sp{\delta}}
	\left(\frac{\log\sp{3}N}{\delta}+\frac{3\log\sp{2}N}{\delta\sp{2}}
	+\frac{6\log\,N}{\delta\sp{3}}+\frac{6}{\delta\sp{4}}\right)\\
&\qquad\qquad +4.005\frac{1+\delta}{N\sp{\delta}}
	\left(\frac{\log\sp{2}N}{\delta}+\frac{2\log\,N}{\delta\sp{2}}
	+\frac{2}{\delta\sp{3}}\right).\\
\endaligned$$
\endproclaim

\linesep

\noindent{\it Proof of Lemma 6.1}.  Using the partial summation formula, one gets
$$\aligned
&\sum\sb{N<n<\infty}\frac{d\sp{2}(n)}{n\sp{2+2\delta}}
      =\int\sb{N}\sp{\infty}\frac{1}{y\sp{2+2\delta}}
           \dd\left(\sum\sb{N<n\le\,y}d\sp{2}(n)\right)\\
&=\left(\frac{1}{y\sp{2+2\delta}}
	\sum_{N<n\le\,y} d\sp{2}(n)\right)\bigg\vert\sb{N}\sp{\infty}
    +(2+2\delta)\int\sb{N}\sp{\infty}\left(\sum\sb{N<n\le\,y} d\sp{2}(n)\right)
            \frac{1}{y\sp{3+2\delta}} \dd y\\
&=(2+2\delta)\int\sb{N}\sp{\infty}\left(\sum_{N<n\leq y} d\sp{2}(n)\right)
            \frac{1}{y\sp{3+2\delta}}\dd y.\\
\endaligned$$
By the Corollary of Lemma 4.2 in [5], one has 
$$\multline
\sum\sb{n\le\,x}d\sp{2}(x)\le 0.102x\log\sp{3}x+1.676x\log\sp{2}x+8.564x\log\,x+23.652x\\
+1.334\sqrt{x}\log\sp{3}x-2.845\sqrt{x}\log\sp{2}x-4.280\sqrt{x}\log\,x-8.501\sqrt{x}\\
+1.334\log\sp{3}x-0.845\log\sp{2}x+2.874\log\,x-0.111\\
\le\,0.103x\log\sp{3}x.\\
\endmultline$$
It follows that 
$$\multline
\sum_{N<n<\infty}\frac{d\sp{2}(n)}{n\sp{2+2\delta}}
\le 0.103(2+2\delta)\int\sb{N}\sp{\infty}\frac{\log\sp{3}y}
	{y\sp{2+2\delta}}\dd\,y.\\
\endmultline\tag{6.1}$$
From this, Lemma 6.1 follows. \qed\par\bigskip

\noindent{\it Proof of Lemma 6.2}. We note that 
$$\sum\sb{N<n}\sum\sb{N<m<n}\frac{d(m)d(n)}{(mn)\sp{1+\delta}}
	\le\left(\sum\sb{N<n}\frac{d(n)}{n\sp{1+\delta}}\right)\sp{2} $$
and by Lemma 5.1 in [5]
$$\sum\sb{n\le\,x}d(n)\le\,x\log\,x+0.155x+4\sqrt{x}\le1.001x\log\,x.$$
Using these, as before, we similarly prove Lemma 6.2. \qed\par\bigskip

\proclaim{Proposition 6.1} For $\log\log\,x\ge18$, one has
$$\sum\sb{m\le\,x}\sum\sb{n<m}\frac{d(m)d(n)}{(mn)\sp{1/2}\log(m/n)}
\le0.066x\log\sp{3}x+4.005x\log\sp{2}x.$$
\endproclaim

\noindent{\it Proof}. Note that $-\log(1-x)>x$ for $0<x<1$. Thus, for $n<m$, 
$$\multline\frac{1}{\log(m/n)}=\left(-\log\bigg(1-{m-n\over m}\bigg)\right)^{-1}\\
<\left( {m-n\over m}\right)^{-1} = 1+{n\over m-n}<1+ {(mn)^{1/2}\over m-n}.\\
\endmultline$$
It follows that 
$$\aligned
&\sum\sb{m\le\,x}\sum\sb{n<m}\frac{d(m)d(n)}{(mn)\sp{1/2}\log(m/n)}\\
&<\sum\sb{m\le\,x}\sum\sb{n<m}\frac{d(m)d(n)}{(mn)\sp{1/2}}
	+\sum\sb{m\le\,x}\sum\sb{n<m}\frac{d(m)d(n)}{m-n}.\\
\endaligned\tag{6.5}$$
For the first sum in (6.5), one sees 
$$\sum\sb{m\le\,x}\sum\sb{n<m} {d(m) d(n) \over (mn)\sp{1/2}}
=\sum\sb{m\le\,x} {d(m)\over m\sp{1/2}} \sum\sb{n\le\,x} {d(n)\over n\sp{1/2}}
\le\left(\sum_{n\leq x} {d(n)\over n^{1/2}}\right)^{2}.$$
Recalling the Corollary of Lemma 5.2 in [5], one has 
$$\sum\sb{n\le x}\frac{d(n)}{\sqrt{n}}
\le2\sqrt{x}\log{x}-1.691\sqrt{x}+2\log\,x+5.846\le2.001\sqrt{x}\log\,x.$$
For the second sum in (6.5), one recalls the Corollary of the main Theorem in [5]. 
Since $\log\log\,x\ge18$, one has 
$$\sum\sb{m\le\,x}\sum\sb{n<m}\frac{d(m)d(n)}{m-n}
\le0.066x\log\sp{3}x.$$

Conclude that one finishes the proof of Proposition 6.1. \qed\par\medskip\bigskip

\noindent{\it Proof of Lemma 6.3}.  Using the partial summation formula for the sum
over $n$, one gets
$$\aligned
& \sum\sb{N<m}\frac{1}{m\sp{1+\delta}}
	\sum\sb{N<n<m}\frac{d(m)d(n)}{(mn)\sp{1/2}\log(m/n)}\\
&\qquad=\int_{N}\sp{\infty}\frac{1}{y\sp{1+\delta}}
        \dd\Big(\sum\sb{N<m\le\,y}\sum\sb{N<n<m}\frac{d(m)d(n)}{(mn)\sp{1/2}
              \log(m/n)}\Big)\\
&\qquad\qquad=\Big(\frac{1}{y\sp{1+\delta}}
	\sum\sb{N<m\le\,y}\sum\sb{N<n<m}\frac{d(m)d(n)}{(mn)\sp{1/2}
	\log(m/n)}\Big)\bigg\vert\sb{N}\sp{\infty}\\
&\qquad\qquad\qquad +(1+\delta)\int\sb{N}\sp{\infty}
	\sum\sb{N<m\le\,y}\sum\sb{N<n<m}\frac{d(m)d(n)}{(mn)\sp{1/2}\log(m/n)}
	\frac{\dd y}{y\sp{2+\delta}}\\
&\qquad\qquad\qquad\qquad=\lim\sb{y\to\infty}\frac{1}{y\sp{1+\delta}}
	\sum\sb{N<m\le\,y}\sum\sb{N<n<m}\frac{d(m)d(n)}{(mn)\sp{1/2}
	\log(m/n)}\\
&\qquad\qquad\qquad\qquad\qquad +(1+\delta)\int\sb{N}\sp{\infty}
	\sum\sb{N<m\le\,y}\sum\sb{N<n<m}\frac{d(m)d(n)}{(mn)\sp{1/2}\log(m/n)}
	\frac{\dd y}{y\sp{2+\delta}}.\\
\endaligned$$
Recalling Proposition 6.1, one sees the first term in the last expression is zero;
and applying (6.4) and (6.3), one obtains 
$$\multline
\sum\sb{N<m}\frac{1}{m\sp{1+\delta}}
     \sum\sb{N<n<m}\frac{d(m)d(n)}{(mn)\sp{1/2}\log(m/n)}\\
\le0.066(1+\delta)\int\sb{N}\sp{\infty}\frac{\log\sp{3}y}{y\sp{1+\delta}}\dd\,y
	+4.005(1+\delta)\int\sb{N}\sp{\infty}\frac{\log\sp{2}y}{y\sp{1+\delta}}\dd\,y\\
\le0.066\frac{1+\delta}{N\sp{\delta}}
	\left(\frac{\log\sp{3}N}{\delta}+\frac{3\log\sp{2}N}{\delta\sp{2}}
	+\frac{6\log\,N}{\delta\sp{3}}+\frac{6}{\delta\sp{4}}\right)\\
+4.005\frac{1+\delta}{N\sp{\delta}}
	\left(\frac{\log\sp{2}N}{\delta}+\frac{2\log\,N}{\delta\sp{2}}
	+\frac{2}{\delta\sp{3}}\right).\\
\endmultline$$
This proves Lemma 6.3. \qed

\head 7. Proofs for Lemma 5.2 and 5.3.\endhead
\noindent{\it Proof of Lemma 5.2.} Recall the definition of $V\sb{A}(s)$ 
from Section 3. Using $(x+y)\sp{2}\le 2(x\sp{2}+y\sp{2})$ for real numbers $x$ 
and $y$, one gets
$$|V\sb{A}(s)|\sp{2}\le 2(|\zeta(s)|\sp{2} |U\sb{A}(s)|\sp{2}+1).\tag{7.1}$$
Recalling the definition of ${\Cal V}\sb{\sigma}(t)$ from Section 5, applying 
Lemma 5.1, and using the same inequality for any $x$ and $y$ again, one acquires  
$${\Cal V}\sb{1/2}(t)\le \left(4C\sp{2}t\sp{2\alpha}\log\sp{2\beta}(t+e)
+4D\sp{2}\right)\int\sb{0}\sp{t}|U\sb{A}(0.5+i\tau)|\sp{2}\dd\tau
+2t.\tag{7.2}$$\par

The integral in the last expression is $$\aligned
&\int\sb{0}\sp{t}|U\sb{A}(0.5+i\tau)|\sp{2}\dd\tau 
=\int\sb{0}\sp{t} U\sb{A}(0.5+i\tau)\overline{U\sb{A}(0.5+i\tau)} 
    \dd\tau\\
&\qquad
=\sum\sb{m=1}\sp{A}\sum\sb{n=1}\sp{A}\frac{\mu(m)\mu(n)}{m\sp{1/2} n\sp{1/2}}
\int\sb{0}\sp{t}\left(\frac{m}{n}\right)\sp{i\tau}\dd\tau.\\
\endaligned$$
Thus, using the inequality
$$\left|\int\sb{0}\sp{t} \left(\frac{m}{n}\right)\sp{i\tau}\dd\tau\right|
\le \frac{2}{\log(m/n)}\quad(m>n),$$ 
we immediately get 
$$\int\sb{0}\sp{t} 
|U\sb{A}(0.5+i\tau)|\sp{2}\dd\tau\le t\sum\sb{n\le A} \frac{1}{n} 
+4\sum\sb{m\le\,A}\sum\sb{n<m}\frac{1}{m\sp{1/2} n\sp{1/2}\log(m/n)}.
\tag{7.3}$$\par

The first term in (7.3) is bounded by $t(\log A+1)$.  For the 
second term on the right side of the last expression, we note
that $x\log x-x+1> 0$ for $x>1$.  It implies that  
$$\frac{1}{\log x}< \frac{x}{x-1} 
= 1+\frac{1}{x-1} < 1+\frac{x\sp{1/2}}{x-1}.$$
We use this for $x=m/n$, getting $$\frac{1}{\log(m/n)} 
< 1+ \frac{n\sp{1/2} m\sp{1/2}}{m-n}.$$  
It follows that 
$$\multline 
\sum\sb{m\le\,A}\sum\sb{n<m}\frac{1}{m\sp{1/2} n\sp{1/2} \log(m/n)}\\
<\sum\sb{m\le\,A}\sum\sb{n<m}\frac{1}{m\sp{1/2} n\sp{1/2}} 
        +\sum\sb{m\le\,A}\sum\sb{n<m}\frac{1}{m-n}\\
\le\left(\sum\sb{n\le\,A}\frac{1}{n\sp{1/2}}\right)\sp{2}
     +\sum\sb{1<n\le A} \Big(1+\log(m-1)\Big)\\
\le\,4A+A\log A\le A(\log A+4).\\
\endmultline$$
Thus, $$\int\sb{0}\sp{t} |U\sb{A}(0.5+i\tau)|\sp{2}\dd\tau
\le t(\log A+1)+4A(\log A+4).\tag{7.4}$$

Conclude that, from (7.2) and (7.4), one shows Lemma 5.2. 

\noindent{\it Proof of Lemma 5.3.}  To estimate ${\Cal V}\sb{1+\delta}(t)$, 
one recalls Lemma 3.1. It follows that 
$${\Cal V}\sb{1+\delta}(t)=\int\sb{0}\sp{t}\left|\sum\sb{A<n}
	\frac{\nu(n)}{n\sp{1+\delta+i\tau}}\right|\sp{2}\dd\,\tau
	=\sum\sb{A<m}\frac{\nu(m)}{m\sp{1+\delta}}
	\sum\sb{A<n}\frac{\nu(n)}{n\sp{1+\delta}}
	\int\sb{0}\sp{t}\left(\frac{m}{n}\right)\sp{i\tau}\dd\,\tau.$$
Similarly to the argument for obtaining (8.3), one deduces 
$${\Cal V}\sb{1+\delta}(t)\le\,t\sum\sb{A<n}\frac{d\sp{2}(n)}{n\sp{2+2\delta}}
    +4\sum\sb{A<m}\sum\sb{A<n<m}\frac{d(m)d(n)}{m\sp{1+\delta} 
	n\sp{1+\delta}\log(m/n)}.\tag{7.5}$$  

For the second sum in the last expression, we observe that 
the function  $f(x)=\log x+x\sp{-1/2}-1>0$ for $x>1$. It 
follows that $$\frac{1}{\log x} < 1+ \frac{1}{x\sp{1/2} \log x}, 
\quad\hbox{for}\ x>1.$$ With $x=m/n$, one sees that 
$$\frac{1}{\log(m/n)}<1+\frac{n\sp{1/2}}{m\sp{1/2}\log(m/n)}.$$
The second term in (7.5) is less than  
$$\aligned &\hskip 2.4 true cm 
\sum\sb{A<m} \sum\sb{A<n<m}  
\frac{d(m) d(n)}{m\sp{1+\delta} n\sp{1+\delta}\log(m/n)}\\
&\le\sum\sb{A<m}\sum\sb{A<n<m}\frac{d(m)d(n)}{m\sp{1+\delta}n\sp{1+\delta}}
+\sum\sb{A<m}\sum\sb{A<n<m}\frac{d(m)d(n)}
	{m\sp{1+\delta}n\sp{\delta}m\sp{1/2}n\sp{1/2}\log(m/n)}\\
&\le\sum\sb{A<m}\sum\sb{A<n<m}\frac{d(m)d(n)}{m\sp{1+\delta}n\sp{1+\delta}}
+\sum\sb{A<m}\frac{1}{m\sp{1+\delta}}\sum\sb{A<n<m}\frac{d(m)d(n)}
	{m\sp{1/2}n\sp{1/2}\log(m/n)}.\\
\endaligned$$ 
Applying Lemma 6.1 for the first term in (7.5) and Lemma 6.2 
for the first term and Lemma 6.3 for the second term in the 
last expression, one proves Lemma 5.2. \bigskip

\head 8. Proofs for Lemma 5.4 and 5.5\endhead\noindent

\noindent{\it Proof of Lemma 5.4}. By the definition of the cosine function 
in complex variable $s=\sigma+it$, one sees that 
$$\aligned\cos\left(\frac{s}{2\tau}\right) 
&=\frac{1}{2}\left(e\sp{-\frac{i\,s}{2\tau}}+e\sp{\frac{i\,s}{2\tau}}\right)\\
&\qquad=\frac{1}{2}\left(e\sp{\frac{t}{2\tau}-i\frac{\sigma}{2\tau}}
	+e\sp{-\frac{t}{2\tau}+i\frac{\sigma}{2\tau}}\right)\\
&\qquad\qquad=\frac{1}{2}e\sp{\frac{t}{2\tau}-i\frac{\sigma}{2\tau}}
	\left(1+e\sp{-\frac{t}{\tau}+i\frac{\sigma}{\tau}}\right).\\
\endaligned$$ 
Since $\frac{1}{2}\le\sigma\le 2$, one has $\displaystyle 0<\frac{\sigma}{\tau}
<\frac{\pi}{4}$ for $\tau\ge e$. It follows that $e\sp{i\frac{\sigma}{\tau}}$ is 
in the first half of the first quadrant so that $1<|1+e\sp{-\frac{t}{\tau}+i\frac{\sigma}{\tau}}|
=\sqrt{1+e\sp{-\frac{t}{\tau}} +2e\sp{-\frac{t}{\tau}}\cos(\sigma/\tau)}<2$. Thus, 
one sees that 
$$\frac{1}{2}e\sp\frac{t}{2\tau}
<\left|\cos\left(\frac{s}{2\tau}\right)\right|
<e\sp\frac{t}{2\tau}.$$ For $t>0$, it is easy to see that
$$\left|\frac{s-1}{s}\right|=\sqrt{1-\frac{2\sigma-1}{\sigma\sp2+t\sp2}}\le 1.$$ 
For $\frac{1}{2}<\sigma\le2$ and $t>14$,
$$\left|\frac{s-1}{s}\right|=\sqrt{1-\frac{2\sigma-1}{\sigma\sp2+t\sp2}}
>\sqrt{1-\frac{2\sigma-1}{\sigma\sp2+14\sp2}}\ge\sqrt{\frac{197}{200}}.$$ 
Conclude that one finishes the proof of Lemma 5.4. 

\noindent{\it Proof of Lemma 5.5.} 
One first note that $H$ is an analytic function so that
$${\Cal H}(\sigma)=2\int\sb{0}\sp{\infty} |H(s)|\sp{2} \dd t.$$
From this equation and the first inequality in Lemma 5.4, one sees 
$${\Cal H}(\sigma)\le\,8\int\sb{0}\sp{\infty} e\sp{-\frac{t}{\tau}} 
|V\sb{A}(s)|\sp{2} \dd t.$$ 
One then uses integration by parts, getting 
$$\aligned &\int\sb{0}\sp{\infty} e\sp{-\frac{t}{\tau}} 
|V\sb{A}(\sigma+it)|\sp{2} \dd t =\int\sb{0}\sp{\infty} 
      e\sp{-\frac{t}{\tau}} 
\dd\left( \int\sb{0}\sp{t}|V\sb{A}(\sigma+iy)|\sp{2}\dd\,y\right)\\
&\qquad=\int\sb{0}\sp{\infty} e\sp{-\frac{t}{\tau}}\dd{\Cal V}\sb{\sigma}(t)
	=e\sp{-\frac{t}{\tau}} {\Cal V}\sb{\sigma}(t)\bigg\vert\sb{0}\sp{\infty} 
+\frac{1}{\tau} \int\sb{0}\sp{\infty} e\sp{-\frac{t}{\tau}} 
    {\Cal V}\sb{\sigma}(t)\dd t. \\     \endaligned$$
Note that ${\Cal V}\sb{\sigma}(0)=0$ by definition. From Lemma 3.3, it is easy
to see that ${\Cal V}\sb{\sigma}(t)\ll\,t\sp{4}$; hence, the first term in the 
last expression is zero. Thus,
$${\Cal H}(\sigma)\le \frac{8}{\tau}\int\sb{0}\sp{\infty} e\sp{-\frac{t}{\tau}} 
{\Cal V}\sb{\sigma} (t)\dd\,t.$$
One then substitutes the variable $t$ by $\tau y$ with the variable $y$ and 
the parameter $\tau$, getting  
$${\Cal H}(\sigma)\le\,8\int\sb{0}\sp{\infty}e\sp{-y}{\Cal V}\sb{\sigma}(\tau y)\dd y.
\tag{8.1}$$\par

To estimate ${\Cal H}\left(\frac{1}{2}\right)$ and ${\Cal H}(1+\delta)$, one uses
Lemma 5.2 and 5.3. One needs to calculate the integrals in the forms of 
$${\Cal J}(a,b):=\int\sb{0}\sp{\infty} e\sp{-y} y\sp{a}\log\sp{b}(y+e)\dd y,$$ 
for the ordered sets $\{a,b\}=\{0,0\}$, $\{1,0\}$, $\left\{\frac{1}{3}, 0\right\}$, 
$\left\{\frac{4}{3}, 0\right\}$, $\left\{\frac{1}{3}, 2\right\}$, 
and $\left\{\frac{4}{3}, 2\right\}$.\par\smallskip

For the first two sets of values for $a$ and $b$, it is easy to see
$${\Cal J}(0,0)=\int\sb{0}\sp{\infty} e\sp{-y}\dd y=1,\quad
\hbox{ and }\quad 
{\Cal J}(1,0)=\int\sb{0}\sp{\infty} e\sp{-y} y\dd y=1,$$
using partial integral formula for the second one. 
one then uses a computation package to get
$${\Cal J}\left(\frac{1}{3}, 0\right)
=\int\sb{0}\sp{\infty} e\sp{-y} y\sp{1/3} \dd y=\Gamma\left(\frac{4}{3}\right)\le 0.893\,,$$
$${\Cal J}\left(\frac{4}{3}, 0\right)
=\int\sb{0}\sp{\infty} e\sp{-y} y\sp{4/3} \dd y=\Gamma\left(\frac{7}{3}\right)\le 1.191\,,$$
$${\Cal J}\left(\frac{1}{3}, 2\right)
=\int\sb{0}\sp{\infty} e\sp{-y} y\sp{1/3}\log\sp{2}(y+e)\dd y
\le 1.220\,,$$ and
$${\Cal J}\left(\frac{4}{3}, 2\right)
=\int\sb{0}\sp{\infty} e\sp{-y} y\sp{4/3}\log\sp{2}(y+e)\dd y
\le 1.881\,.$$

Note that $\tau y+e\le\tau(y+e)$ since $\tau\ge e$; so that 
$\log(\tau y+e)\le\log\tau+\log(y+e)$. One then has 
$$\log\sp{2}(\tau y+e)\le 2\left(\log\sp{2}(\tau)+\log\sp{2}(y+e) \right),$$ 
since $(x+y)\sp{2}\le 2(x\sp{2}+y\sp{2})$ is valid for any real numbers $x$ and $y$. 
Recalling~Lemma~5.2, one obtains that 
$$\aligned \int\sb{0}\sp{\infty} e\sp{-y} {\Cal V}\sb{1/2}(\tau y)\dd y
&\le 2 D\sb{1}\tau\sp{4/3}\log\sp{2}\tau\; {\Cal J}\left(\frac{4}{3},0\right)
+2 D\sb{1}\tau\sp{4/3}\; {\Cal J}\left(\frac{4}{3}, 2\right)\\
&\hskip 1 true cm
+2 D\sb{2}\tau\sp{1/3}\log\sp{2}\tau\; {\Cal J}\left(\frac{1}{3},0\right)
	+2 D\sb{2}\tau\sp{1/3}{\Cal J}\left(\frac{1}{3},2\right)\\
&\hskip 2 true cm 
  +D\sb{3}\tau{\Cal J}(1,0)+D\sb{4}{\Cal J}(0,0).\\
\endaligned$$
Then, recalling (8.1), one acquires 
$${\Cal H}\left(\frac{1}{2}\right)
\le a\sb{1} \tau\sp{4/3}\log\sp{2}\tau +a\sb{2}\tau\sp{4/3}
+a\sb{3}\tau\sp{1/3}\log\sp{2}\tau +a\sb{4}\tau\sp{1/3}
+a\sb{5}\tau+a\sb{6},\tag{8.2}$$
with $$\multline
a\sb{1}:=16 D\sb{1} {\Cal J}\left(\frac{4}{3},0\right) \le 14.288 D\sb{1},\quad
	a\sb{2}:=16 D\sb{1} {\Cal J}\left(\frac{4}{3},2\right) \le 19.056 D\sb{1},\\
a\sb{3}:=16 D\sb{2} {\Cal J}\left(\frac{1}{3},0\right) \le 19.520 D\sb{2},\quad 
	a\sb{4}:=16 D\sb{2} {\Cal J}\left(\frac{1}{3},2\right) \le 30.096 D\sb{2},\\
\quad\quad a\sb{5}:=8 D\sb{3} {\Cal J}\left(1,0\right) =8 D\sb{3},\quad\quad 
	a\sb{6}:=8 D\sb{4} {\Cal J}\left(0,0\right) =8 D\sb{4}.\\
\endmultline$$
Similarly, but recalling Lemma 5.3 and (8.1), one has 
$${\Cal H}(1+\delta)\le b\sb{1}\tau +b\sb{2},\tag{8.3}$$ 
with $b\sb{1}=8 D\sb{5}{\Cal J}(1,0) =8 D\sb{5}$ 
and $b\sb{2}=8 D\sb{6}{\Cal J}(1,0) =8 D\sb{6}$.\par\bigskip

Actually, the ``constants'' $a\sb{j}$ for $j=1$, $\ldots$, $6$ and $b\sb{j}$ 
for $j=1$, $2$, are not absolute constants; they depend on the choice of $A$ 
subject to $A\ge 16$ as well as our choice of the parameter $\tau$. 
The kink is that we are going to choose suitable $A$ and $\tau$.\par

Note that $\big(1+\frac{1}{T\sb{0}}\big)T-T=\frac{T}{T\sb{0}}\ge\,1$. One may choose $A$ to be 
an integer in $T\le\,A\le\,\big(1+\frac{1}{T\sb{0}}\big) T$. Let $\kappa$ be a constant such that 
$\kappa\ge\,\frac{e}{T\sb{0}}$ and $\tau=\kappa\,T$. Then $\tau\ge\,e$. Also, let $\omega>0$
and $\delta=\frac{\omega}{T}$.\par

For brevity, denote $A(T) =\log\,T +\log\Big(1+\frac{1}{T\sb{0}}\Big)$. We have $A(T) +Z 
<1.000,000,001\log T$ for any $Z=1$ or $4$. Also, we assume that $\kappa$ is not so large 
so that $\log T +\kappa \le 1.000,001\log T$. It is now straightforward to conclude Lemma 5.5.

\head 9. Landau's Approximate Formula\endhead\noindent 
In this section, we give an explicit form of Landau's approximate formula as stated 
in Lemma 9.1.\par

Let $T\ge0$ and $u>0$. Suppose there are $n$ zeros $\beta\sb{1}+iz\sb{1}$, $\beta\sb{2}+iz\sb{2}$, 
$\ldots$, $\beta\sb{n}+iz\sb{n}$ of $\zeta(s)$ in $T-u\le\Im(s)\le T+u$ such that 
$z\sb{0}=T-u\le z\sb{1}<z\sb{2}<\ldots<z\sb{n}\le T+u=z\sb{n+1}$.
Let $1\le\,j\le\,n+1$ be such that $z\sb{j}-z\sb{j-1}\ge\,z\sb{i}-z\sb{i-1}$ for 
every other $1\le\,i\le\,n+1$. There may be more than one such a $j$. Fix one such $j$
and let $T\sb{u}=\frac{z\sb{j-1}+z\sb{j}}{2}$. For convenience, $T\sb{u}$ is called 
the associate of $T$ with respect to $u$.\par 

\proclaim{Lemma 9.1} Let $x\ge x\sb{0}$ and $T\ge\exp(\exp(18))$. Suppose that $T\sb{u}$ 
is the associate of $T$ with respect to $u=1.155$. Then,
$$\psi(x)=x-\sum\sb{|\Im(\rho)|\le T\sb{u}} \frac{x\sp{\rho}}{\rho}+E(x),$$
where $$|E(x)|\le 5.26\,{x\log\sp{2}x \over T}
	+33.488 \frac{x\log\sp{2}T}{T\log x} +3\frac{\log\sp{2}T}{x}. $$
\endproclaim

\linesep

\proclaim{Proposition 9.1}
Let $t\ge\,0$ and $\beta\sb{n}+i\gamma\sb{n}$, $n=1$, $2$, $\ldots$ be all non-trivial 
zeros of the Riemann zeta-function.  Then
$$\sum\sb{n=1}\sp{\infty}\frac{1}{4+(t-\gamma\sb{n})\sp{2}}
\le\frac{1}{4}\log\left(t\sp{2}+4\right)+1.483\,.$$ 
\endproclaim

\noindent{\it Proof}. Recall the following formula, see [8]. That is,  
$$-{\zeta'(s)\over \zeta(s)}=-{1\over s-1} +\sum_{n=1}^{\infty} \left({1\over s-\rho_{n}}
+{1\over \rho_{n}} \right) +\sum_{n=1}^{\infty} \left({1\over s+2n}
-{1\over 2n}\right) +B_{0},\tag{9.1}$$ where $\{ \rho_{n}: n=1,2, \ldots\}$ is the set
of all non-trivial zeros of the Riemann zeta-function and $B\sb{0}=\log(2\pi)-1$. 
Using this equation with $s=2+it$, Proposition 9.1 follows.  

\proclaim{Proposition 9.2}
Let $\gamma_{n}$ be defined in Lemma \hbox{\rm 9.1}. For $t\ge\,0$ and $0<u$, 
one has\par\medskip 
\noindent \hbox{\rm (a)} The number of zeros of $\zeta(s)$ such that
$|t-\gamma_{n}|\leq u$ is less than 
$$(4+u^{2})\left(\frac{1}{4}\log\left(t\sp{2}+4\right)+1.483\right);$$
\par\smallskip
\noindent\hbox{\rm (b)} $\displaystyle{\sum_{|t-\gamma_{n}|>u}
\frac{1}{(t-\gamma_{n})\sp{2}}\le\left(1+\frac{4}{u\sp{2}}\right)
\left(\frac{1}{4}\log\left(t\sp{2}+4\right)+1.483\right)}$.
\endproclaim

\noindent{\it Proof}. Note that 
$$1\le{4+u\sp{2}\over 4+(t-\gamma_{n})^{2}},
\quad \hbox{ if $|t-\gamma_{n}|\leq u$,} $$ therefore,
$$\sum_{|t-\gamma_{n}|\leq u} 1\le (4+u^{2})
\sum_{n=1}^{\infty} {1\over 4+(t-\gamma_{n})^{2}}.$$
Applying Proposition 9.1, one proves (a) in Proposition 9.2. One shows (b) 
in the proposition similarly, but note that
$${1\over (t-\gamma_{n})^{2}}\leq \left(1+{4\over u^{2}}\right)
{1\over 4+(t-\gamma_{n})^{2}}, \quad \hbox{ if $|t-\gamma_{n}|> u$, }$$
so that
$$\sum\sb{|t-\gamma_{n}|>u} \frac{1}{(t-\gamma)\sp{2}}
\le\left(1+\frac{4}{u^{2}}\right)
\sum\sb{n=1}\sp{\infty} \frac{1}{4+(t-\gamma\sb{n})\sp{2}}.\qed$$
\par\bigskip

\proclaim{Proposition 9.3} Let $-1\le \sigma\le 2$, $t>0$, and $u>0$.
Then $$\multline
\left|\frac{\zeta\sp{\prime}(\sigma\pm it)}{\zeta(\sigma\pm it)}\right|
\le\sum\sb{|t-\gamma\sb{n}|\le\,u}\left(\frac{1}{2+it-\rho\sb{n}}
	-\frac{1}{s-\rho\sb{n}}\right)\\
+\frac{3}{2}\left(1+\frac{4}{u\sp{2}}\right)
	\left(\frac{1}{4}\log\Big(t\sp{2}+4\Big)+1.483\right)+\frac{3}{t\sp{2}}+1.284.\\
\endmultline$$
\endproclaim

\noindent{\it Proof}. From (9.1), one has the following equation.
$$\eqalign{&-{\zeta'(s)\over\zeta(s)} 
	=\frac{1}{2+it-1}-\frac{1}{s-1}
	+\sum\sb{n=1}\sp{\infty}\left(\frac{1}{s-\rho_{n}}
	-\frac{1}{2+it-\rho_{n}}\right)\cr
&\hskip 2 true cm +\sum_{n=1}^{\infty} 
	\left(\frac{1}{s+2n}-\frac{1}{2+it+2n}\right)
	-{\zeta'(2+it)\over\zeta(2+it)}. \cr} $$
Using this equation and Proposition 9.2(b), Proposition 9.3 follows.  

\proclaim{Proposition 9.4}Let $-1\le\sigma\le 2$ and $T>exp(exp(18))$. Suppose 
$T\sb{u}$ is the associate of $T$ with respect to $u=1.155$. Then\par
\noindent \hbox{\rm (a)}
$$\left|\frac{\zeta\sp{\prime}(\sigma\pm\,iT\sb{u})}{\zeta(\sigma\pm\,iT\sb{u})}\right|
\le6.159\log\sp{2}T+2.999\log T+1.285;$$
\noindent \hbox{\rm (b)} For $12<t\le\,T\sb{u}$ 
$$\left|\frac{\zeta\sp{\prime}(-1\pm\,it)}{\zeta(-1\pm\,it)}\right|
\le2.999\log\,t+10.241;$$  
\par\smallskip
\noindent
and\par 
\noindent \hbox{\rm (c)} For $0\le\,t\le12$ 
$$\left|\frac{\zeta\sp{\prime}(-1\pm\,it)}{\zeta(-1\pm\,it)}\right|\le\,19.172.$$  
\endproclaim

\noindent{\it Proof}. Note that 
$$\left|\frac{\zeta\sp{\prime}(\sigma-it)}{\zeta(\sigma-it)}\right|=
\left|\frac{\zeta\sp{\prime}(\sigma+it)}{\zeta(\sigma+it)}\right|.$$
One only needs to consider the case with the plus sign for each case. \par

Recalling Proposition 9.3, one only needs to estimate the sum
$$\sum\sb{|t-\gamma\sb{n}|\le\,u}\left(\frac{1}{s-\rho\sb{n}}
	-\frac{1}{2+it-\rho\sb{n}}\right).\tag{9.2}$$
Note that 
$$\frac{1}{|s-\rho\sb{n}|}
=\frac{1}{|\sigma-\beta\sb{n}+i(t-\gamma\sb{n})|}
=\frac{1}{\sqrt{(\sigma-\beta\sb{n})\sp{2}+(t-\gamma\sb{n})\sp{2}}}
\le\frac{1}{|t-\gamma\sb{n}|}.\tag{9.3}$$
Similarly, 
$$\frac{1}{|2+it-\rho\sb{n}|}\le\frac{1}{|t-\gamma\sb{n}|}.\tag{9.4}$$
Thus, 
$$\left|\sum\sb{|t-\gamma\sb{n}|\le\,u}\left(\frac{1}{\sigma+it-\rho\sb{n}}
	-\frac{1}{2+it-\rho\sb{n}}\right)\right|
\le2\sum\sb{|t-\gamma\sb{n}|\le\,u}\frac{1}{|t-\gamma\sb{n}|}.\tag{9.5}$$
Recalling (a) in Proposition 9.2, one sees that there are at most 
$$(4+u\sp{2})(\log(t\sp{2}+4)/4+1.483)$$ terms in (9.5). By the setting of 
$T\sb{u}$, one has
$$|T\sb{u}-\gamma\sb{n}|\ge\frac{2u}{(4+u\sp{2})(\log(T\sp{2}+4)/4+1.483)+1}$$ 
for every $\gamma\sb{n}$, $n=1$, $2$, $\ldots$ (It is $T$ instead of $T\sb{u}$ on
the right side of the last expression). Or, each summand in the last expression 
in (9.5) is less than $$\frac{(4+u\sp{2})(\log(T\sp{2}+4)/4+1.483)+1}{2u}.$$
It follows that 
$$\multline\left|\sum\sb{|T\sb{u}-\gamma\sb{n}|\le\,u}
	\left(\frac{1}{\sigma+iT\sb{u}-\rho\sb{n}}
	-\frac{1}{2+iT\sb{u}-\rho\sb{n}}\right)\right|\\
\le2(4+u\sp{2})(\log(T\sb{u}\sp{2}+4)/4+1.483)
	\frac{(4+u\sp{2})(\log(T\sp{2}+4)/4+1.483)+1}{2u}\\
\le(4+u\sp{2})(\log((T+u)\sp{2}+4)/4+1.483)
	\frac{(4+u\sp{2})(\log(T\sp{2}+4)/4+1.483)+1}{u}.\\
\endmultline$$
To sub-optimize the factor $\frac{(4+u\sp{2})\sp{2}}{u}$, one let $u=1.155$. 
Also, remark such as  
$$\multline
\qquad\qquad\log(T\sp{2}+4)=2\log\,T+\log\left(1+4e\sp{e\sp{-36}}\right)\le2.0000001\log\,T,\\
\hbox{and}\\
\log((T+u)\sp{2}+4)=2\log\,T+\log\left(\big(1+1.155e\sp{e\sp{-18}}\big)\sp{2}
	+4e\sp{e\sp{-36}}\right)\le2.0000001\log\,T.\\ 
\endmultline$$
Summarizing with the result in Proposition 9.3, one proves (a). \par

One proves (b) and (c) similarly. For (b), replacing (9.3) and (9.4) by
$$\frac{1}{|-1+it-\rho\sb{n}|}\le\frac{1}{|-1-\beta\sb{n}|}\le 1,\quad\hbox{and}\quad
\frac{1}{|2+it-\rho\sb{n}|}\le\frac{1}{|2-\beta\sb{n}|}\le 1,$$
and noting that 
$$\log(t\sp{2}+4)=2\log\,t+\log(1+4/12\sp{2})\le2\log\,t+0.028.$$
For (c), one replaces the upper bound in (9.8) by $\frac{3}{2}$, recalling 
(13) from [7]. Also, note that the terms of sum in (9.2) is zero and 
$$\log(t\sp{2}+4)\le\log(12\sp{2}+4)\le4.998.\qed$$

We also need Landau's Approximate Formula in the following form, see Lemma~4 from [7]. \par

\proclaim {Proposition 9.5}. Let $x>2981$ and $T\ge\exp(\exp(18))$. Suppose that 
$T\sb{u}$ is the associate of $T$ with respect to $u=1.155$.
Then $$\psi(x)={1\over 2\pi i}
\int\sb{1+\frac{1}{\log x}-iT\sb{u}}\sp{1+\frac{1}{\log x}+iT\sb{u}} 
\left(-\frac{\zeta'(s)}{\zeta(s)} \right)\frac{x\sp{s}}{s}\dd s+E\sb{0}(x),$$
where $$|E\sb{0}(x)|\le 5.25?\,{x\log\sp{2}x \over T}+12.64?\,{x\log x \over T}
+\log x.$$
Especially, if $x\ge e\sp{e\sp{15}}$, then 
$$|E\sb{0}(x)|< 5.26\,{x\log\sp{2}x \over T} -\log(2\pi x).$$
\endproclaim\bigskip

\noindent{\it Proof of Lemma 9.1}. One applies Cauchy's Residue Theorem on the function 
$-\frac{\zeta\sp{\prime}(s)}{\zeta(s)}\frac{x\sp{s}}{s}$. Utilizing (9.1), 
we see that the residue of $-\frac{\zeta\sp{\prime}(s)}{\zeta(s)}\frac{x\sp{s}}{s}$ 
at $s=1$ is $x$, those at $s=\rho\sb{n}$'s are $-\frac{x\sp{\rho}}{\rho}$'s and 
that at $s=0$ is $-\frac{\zeta\sp{\prime}(s)}{\zeta(0)}$. We let $r=1+\frac{1}{\log x}
<1.01$ as in section 4 in [7] and $l=-1$ and apply Cauchy's Residue Theorem on the 
rectangle bounded by $s=l$, $s=-iT\sb{u}$, $s=r$, and $s=iT\sb{u}$, getting
$$\aligned &\hskip 2true cm\frac{1}{2\pi i} \int\sb{r-iT\sb{u}}\sp{r+iT\sb{u}}
\left(-\frac{\zeta\sp{\prime}(s)}{\zeta(s)} \frac{x\sp{s}}{s}\right) \dd s\\
&=x-\sum\sb{|\Im(\rho)|\le\,T\sb{u}} 
	\frac{x\sp{\rho}}{\rho} +\frac{\zeta\sp{\prime}(0)}
	{\zeta(0)}-\frac{1}{2\pi i}\int\sb{L\sb{l}+L\sb{u}+L\sb{b}}
\left(-\frac{\zeta\sp{\prime}(s)}{\zeta(s)} \frac{x\sp{s}}{s}\right) \dd s,\\ 
\endaligned\tag{9.6}$$
where $L\sb{l}$ is the left, $L\sb{u}$ is the top, and $L\sb{b}$ is the bottom sides of 
the rectangle. For the third term on the right side of (9.6), one has
$$\frac{\zeta'(0)}{\zeta(0)}=B\sb{0}+1=\log(2\pi).$$\par

For the integral along with $L\sb{l}$, one uses $|-1+it|\ge 1$ for $|t|\le12$ 
and $|-1+it|\ge t$ otherwise to get  
$$\multline
\left|\int\sb{L\sb{l}}\left(-\frac{\zeta\sp{\prime}(s)}
	{\zeta(s)}\frac{x\sp{s}}{s}\right)\dd s\right|
\le\int\sb{-T\sb{u}}\sp{T\sb{u}} \left|\left(-\frac{\zeta\sp{\prime}
	(-1+it)}{\zeta(-1+it)}\right)\frac{x\sp{-1}}{-1+it}\right|\dd t\\
\le\frac{2}{x}\int\sb{12}\sp{T\sb{u}}\frac{2.999\log t+10.241}{t}\dd t
	+\frac{38.344}{x}\int\sb{0}\sp{12}\dd t\\
=\frac{2.999\log\sp{2}T\sb{u}+20.482\log T\sb{u}+460.128-2.999\log\sp{2}12
	-20.482\log\,12}{x} \\
\le\frac{2.999\log\sp{2}(T+1.155)+20.482\log(T+1.155)+390.715}{x} 
	<3\frac{\log\sp{2}T}{x}. \\  
\endmultline$$
For the integral along with $L\sb{u}$ and $L\sb{b}$, one has 
$$\multline
\left|\int\sb{L\sb{l}}
\left(-\frac{\zeta\sp{\prime}(s)}{\zeta(s)}\frac{x\sp{s}}{s}\right)\dd s
	+\int\sb{L\sb{b}}\left(-\frac{\zeta\sp{\prime}(s)}
	{\zeta(s)}\frac{x\sp{s}}{s}\right)\dd s\right|\\
\le2\int\sb{-1}\sp{1+\frac{1}{\log\,x}}\left|
	\left(-\frac{\zeta\sp{\prime}(\sigma+iT\sb{u})}
	{\zeta(\sigma+iT\sb{u})}\right)\frac{x\sp{\sigma+iT\sb{u}}}
	{\sigma+iT\sb{u}}\right|\dd\sigma\\
\le2\frac{6.159\log\sp{2}T+2.999\log T+1.285}{T\sb{u}} 
	\int\sb{-1}\sp{1+\frac{1}{\log\,x}} x\sp{\sigma}\dd\sigma\\
=\frac{(12.318\log\sp{2}T+5.998\log T+2.570)(ex-1/x)}{T\sb{u}\log\,x}\\
\le\frac{ex(12.318\log\sp{2}T+5.998\log T+2.570)}{(T-u)\log x}\\ 
\le\frac{ex(12.319\log\sp{2}T+5.999\log T+2.571)}{T\log x}
	< 33.488 \frac{x\log\sp{2}T}{T\log x}. \\ 
\endmultline$$
Conclude that one proves Lemma 9.1.

\enddocument